\documentclass[leqno,11pt]{article}
\usepackage{amssymb,amsmath}
\renewcommand\Box{{\vrule height5pt width 5pt depth0pt}}

\newcommand{\Z}{{\mathbb{Z}}}
\newcommand{\Q}{{\mathbb{Q}}}
\newcommand{\N}{{\mathbb{N}}}
\newcommand{\C}{{\mathbb{C}}}
\newcommand{\fp}{{\mathfrak{p}}}
\newcommand{\fq}{{\mathfrak{q}}}
\newcommand{\Sym}{{\mathfrak{S}}}
\newcommand{\cB}{{\mathcal{B}}}
\newcommand{\bH}{{\mathbf{H}}}
\newcommand{\Cl}{{\mbox{\rm Cl}}}
\newcommand{\Irr}{{\mbox{\rm Irr}}}

\renewcommand{\geq}{\geqslant}
\renewcommand{\leq}{\leqslant}

\newtheorem{thm}{Theorem}[section]
\newtheorem{lem}[thm]{Lemma}

\newtheorem{prop}[thm]{Proposition}

\newtheorem{defn}[thm]{Definition}
\newtheorem{rema}[thm]{Remark}

\newenvironment{abschnitt}{\refstepcounter{thm} 
\medskip \noindent ({\bf \thethm})}{\par\medskip}

\newenvironment{proof}{\noindent {\it Proof.}}{\hfill $\Box$ \smallskip}

\begin{document}

\title{{\normalsize\it To Maria Pop on her $60$th birthday}\\[5mm]
On the number of simple modules of Iwahori--Hecke algebras of finite
Weyl groups}

\author{Meinolf Geck}
\date{}
\maketitle

\pagestyle{myheadings}
\markboth{}{Geck: Number of simple modules}

\begin{abstract} 
Let $H_k(W,q)$ be the Iwahori--Hecke algebra associated with a finite
Weyl group $W$, where $k$ is a field and $0 \neq q \in k$. Assume that the
characteristic of $k$ is not ``bad'' for $W$ and let $e$ be the smallest
$i \geq 2$ such that $1+q+q^2+\cdots +q^{i-1}=0$. We show that the
number of simple $H_k(H,q)$-modules is ``generic'', i.e., it only depends 
on~$e$. The proof uses some computations in the {\sf CHEVIE} package of 
{\sf GAP} and known results due to Dipper--James, Ariki--Mathas, Rouquier 
and the author.
\end{abstract}

\renewcommand{\thefootnote}{}

\section{Introduction}
Let $\Gamma$ be one of the graphs in Table~1.\footnote{1991 {\it Mathematics 
Subject Classification.} Primary 20C20. Secondary 20G05.} Let $S$ be the 
set of nodes of $\Gamma$. For $s,t \in S$, $s\neq t$, we define an integer 
$m(s,t)$ as follows. If $s,t$ are not joined in the graph, then $m(s,t)=2$; 
if $s,t$ are joined by an unlabelled edge, then $m(s,t)=3$. Finally, if 
$s,t$ are joined by an edge with label $m$, we set $m(s,t)=m$. We then 
define a group $W=W(\Gamma)$ by the following presentation:
\[ W=\langle s \in S \mid \mbox{$s^2=1$ for $s \in S$ and $(st)^{m(s,t)}=1$
for $s,t \in S, s \neq t$}\rangle.\]
It is known that the groups defined in this way are all finite and that 
they are precisely the Weyl groups arising in the theory of complex simple 
Lie algebras or linear algebraic groups; see \cite{Ca}. For example, if 
$\Gamma$ is the graph $A_{n}$, then $W$ is isomorphic to the symmetric 
group $\Sym_{n+1}$, where the generator $s_i \in S$ attached to the node 
labelled by $i$ corresponds to the basic transposition $(i,i+1)$ 
for $1 \leq i \leq n$.

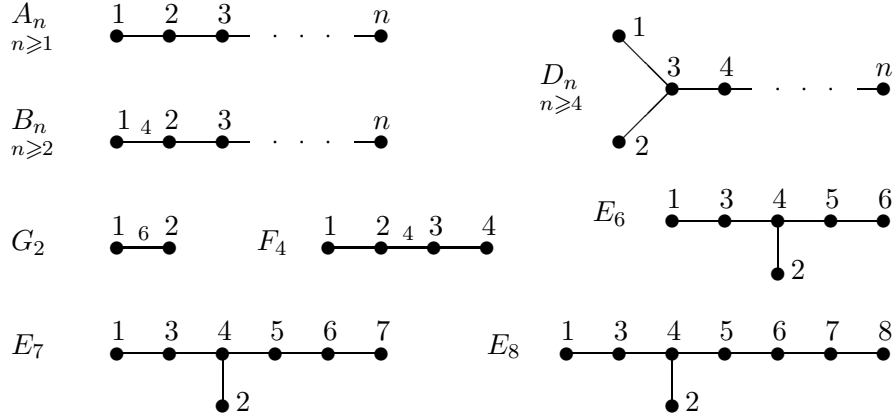
\begin{table}[htbp] \caption{Graphs of finite Weyl groups} 
\begin{center}
\unitlength=1pt
\makeatletter
\vbox{\begin{picture}(342,165)
\put( 10, 25){$E_7$}
\put( 50, 25){\@dot{5}}
\put( 48, 30){$1$}
\put( 50, 25){\line(1,0){20}}
\put( 70, 25){\@dot{5}}
\put( 68, 30){$3$}
\put( 70, 25){\line(1,0){20}}
\put( 90, 25){\@dot{5}}
\put( 88, 30){$4$}
\put( 90, 25){\line(0,-1){20}}
\put( 90,  5){\@dot{5}}
\put( 95,  3){$2$}
\put( 90, 25){\line(1,0){20}}
\put(110, 25){\@dot{5}}
\put(108, 30){$5$}
\put(110, 25){\line(1,0){20}}
\put(130, 25){\@dot{5}}
\put(128, 30){$6$}
\put(130, 25){\line(1,0){20}}
\put(150, 25){\@dot{5}}
\put(148, 30){$7$}

\put(190, 25){$E_8$}
\put(220, 25){\@dot{5}}
\put(218, 30){$1$}
\put(220, 25){\line(1,0){20}}
\put(240, 25){\@dot{5}}
\put(238, 30){$3$}
\put(240, 25){\line(1,0){20}}
\put(260, 25){\@dot{5}}
\put(258, 30){$4$}
\put(260, 25){\line(0,-1){20}}
\put(260,  5){\@dot{5}}
\put(265,  3){$2$}
\put(260, 25){\line(1,0){20}}
\put(280, 25){\@dot{5}}
\put(278, 30){$5$}
\put(280, 25){\line(1,0){20}}
\put(300, 25){\@dot{5}}
\put(298, 30){$6$}
\put(300, 25){\line(1,0){20}}
\put(320, 25){\@dot{5}}
\put(318, 30){$7$}
\put(320, 25){\line(1,0){20}}
\put(340, 25){\@dot{5}}
\put(338, 30){$8$}

\put( 10, 64){$G_2$}
\put( 50, 65){\@dot{5}}
\put( 48, 70){$1$}
\put( 58, 68){$\scriptstyle{6}$}
\put( 50, 65){\line(1,0){20}}
\put( 70, 65){\@dot{5}}
\put( 68, 70){$2$}

\put(103, 64){$F_4$}
\put(130, 65){\@dot{5}}
\put(128, 70){$1$}
\put(130, 65){\line(1,0){20}}
\put(150, 65){\@dot{5}}
\put(148, 70){$2$}
\put(150, 65){\line(1,0){20}}
\put(158, 68){$\scriptstyle{4}$}
\put(170, 65){\@dot{5}}
\put(168, 70){$3$}
\put(170, 65){\line(1,0){20}}
\put(190, 65){\@dot{5}}
\put(188, 70){$4$}

\put(230, 75){$E_6$}
\put(260, 75){\@dot{5}}
\put(258, 80){$1$}
\put(260, 75){\line(1,0){20}}
\put(280, 75){\@dot{5}}
\put(278, 80){$3$}
\put(280, 75){\line(1,0){20}}
\put(300, 75){\@dot{5}}
\put(298, 80){$4$}
\put(300, 75){\line(0,-1){20}}
\put(300, 55){\@dot{5}}
\put(305, 53){$2$}
\put(300, 75){\line(1,0){20}}
\put(320, 75){\@dot{5}}
\put(318, 80){$5$}
\put(320, 75){\line(1,0){20}}
\put(340, 75){\@dot{5}}
\put(338, 80){$6$}

\put( 10,110){$B_n$}
\put( 10,100){$\scriptstyle{n \geq 2}$}
\put( 50,105){\@dot{5}}
\put( 50,110){$1$}
\put( 50,105){\line(1,0){20}}
\put( 59,108){$\scriptstyle{4}$}
\put( 70,105){\@dot{5}}
\put( 68,110){$2$}
\put( 70,105){\line(1,0){30}}
\put( 90,105){\@dot{5}}
\put( 88,110){$3$}
\put(110,105){\@dot{1}}
\put(120,105){\@dot{1}}
\put(130,105){\@dot{1}}
\put(140,105){\line(1,0){10}}
\put(150,105){\@dot{5}}
\put(147,110){$n$}

\put( 10,150){$A_n$}
\put( 10,140){$\scriptstyle{n \geq 1}$}
\put( 50,145){\@dot{5}}
\put( 48,150){$1$}
\put( 50,145){\line(1,0){20}}
\put( 70,145){\@dot{5}}
\put( 68,150){$2$}
\put( 70,145){\line(1,0){30}}
\put( 90,145){\@dot{5}}
\put( 88,150){$3$}
\put(110,145){\@dot{1}}
\put(120,145){\@dot{1}}
\put(130,145){\@dot{1}}
\put(140,145){\line(1,0){10}}
\put(150,145){\@dot{5}}
\put(147,150){$n$}

\put(210,127){$D_n$}
\put(210,117){$\scriptstyle{n \geq 4}$}
\put(240,145){\@dot{5}}
\put(245,145){$1$}
\put(240,105){\@dot{5}}
\put(246,100){$2$}
\put(240,145){\line(1,-1){21}}
\put(240,105){\line(1,1){21}}
\put(260,125){\@dot{5}}
\put(258,130){$3$}
\put(260,125){\line(1,0){30}}
\put(280,125){\@dot{5}}
\put(278,130){$4$}
\put(300,125){\@dot{1}}
\put(310,125){\@dot{1}}
\put(320,125){\@dot{1}}
\put(330,125){\line(1,0){10}}
\put(340,125){\@dot{5}}
\put(337,130){$n$}
\end{picture}}
\makeatother
\end{center}
\end{table}

Now let $k$ be any field and $q \in k$ be any non-zero element which has
a square root in~$k$. Then we define an associative $k$-algebra $H=H_k(W,q)$ 
(with identity $1_H$) by a presentation with 
\[ \begin{array}{ll} \mbox{generators:} & T_w,\, w \in W, 
\\ \mbox{relations:} & T_s^2=q\, 1_H+(q-1)T_s \mbox{ for all $s \in S$,}
\\ & T_wT_{w'}=T_{ww'} \mbox{ for $w,w' \in W$ with $l(ww')=l(w)+l(w')$}.
\end{array}\]
Here, the length function $l \colon W \rightarrow {\N}_0$ is defined as 
follows. Any element $w \in W$ can be written in the form $w=s_1 \cdots 
s_m$ with $s_1,\ldots,s_m \in S$; if $m$ is as small as possible, we set 
$l(w):=m$. It is known that $\{T_w \mid w \in W\}$ is in fact a $k$-basis 
of $H$ and that we have the following multiplication rules. Let $s \in S$ 
and $w\in W$. Then we have $l(sw)=l(w) \pm 1$ and 
\[ T_sT_w=\left\{\begin{array}{cl} T_{sw} & \quad \mbox{if $l(sw)=l(w)+1$},\\
qT_{sw}+(q-1)T_w & \quad \mbox{if $l(sw)=l(w)-1$};\end{array}\right.\]
see \cite[Chapter~4]{ourbuch}. These algebras play an important role, for
example, in the representation theory of finite groups of Lie type (see 
\cite{Ca} and \cite{dghm}) or in the theory of knots and links (see 
\cite[Chapter~4]{ourbuch}).

Now the structure of $H$ is relatively well-understood in the case where
$H_k(W,q)$ is split semisimple. Then Lusztig (see \cite{Lu2} and the
references there) has constructed a canonical isomorphism from $H$ onto 
the so-called asymptotic algebra which itself is isomorphic to the group 
algebra of $W$ over $k$ (if the characteristic of $k$ does not divide the 
order of $W$). The case where $H_k(W,q)$ is not semisimple is much more 
difficult and far from being solved. Let $\Irr(H_k(W,q))$ be the set of 
simple $H_k(W,q)$-modules modulo isomorphism. The purpose of this paper is 
to establish the following result.

\begin{thm} \label{thm1} Let $\ell \geq 0$ be the characteristic of $k$. 
Assume that 
\[\begin{array}{ll}
\ell \neq 2 & \quad \mbox{if $W$ is of type $B_n$ or $D_n$}; \\
\ell \neq 2,3 & \quad \mbox{if $W$ is of type $G_2$, $F_4$, $E_6$ or $E_6$};\\
\ell \neq 2,3,5 & \quad \mbox{if $W$ is of type $E_8$}.\end{array}\]
Let $e$ be the smallest $i \geq 2$ such that $1+q+q^2+\cdots +q^{i-1}=0$.
(We set $e=\infty$ if no such $i$ exists.) Then we have 
\[ |\Irr(H_k(W,q))|=|\Irr(H_{\C}(W,\zeta_e))|,\]
where $\zeta_e \in \C$ is a primitive $e$th root of unity. (We set
$\zeta_\infty=1$.) In particular, $|\Irr(H_k(W,q))|$ only depends on $e$, 
but not on the particular choice of $k$ or~$q$.
\end{thm}

The above result is already known in the following cases. 
\begin{itemize}
\item[$A_n$:] In this case, $|\Irr(H_k(W,q))|$ is the number of $e$-regular 
partitions of $n+1$; see Dipper--James \cite{DJ}. A partition of $n+1$,
written in exponential form $(1^{n_1},2^{n_2},\ldots)$, is called $e$-regular 
if $n_i<e$ for all~$i$.
\item[$B_n$:] In this case, $|\Irr(H_k(W,q))|$ is the number of so-called
Kleshchev bipartitions of $n$; see Ariki--Mathas \cite{AM} for the proof 
and the precise definition.
\item[$D_n$:] In this case, the problem of computing $|\Irr(H_k(W,q))|$ can 
be reduced to the analogous problem for type $B_n$ (but one has to consider
an algebra of type $B_n$ with unequal parameters); see 
Geck~\cite{my00}.
\end{itemize}
For the exceptional types $G_2$, $F_4$ and $E_6$, the numbers 
$|\Irr(H_k(W,q))|$ are explicitly known by work of Geck and Lux without any 
restriction on the characteristic, but assuming that $\ell=0$ or that $q$ 
lies in the prime field of $k$; see \cite[\S 3.4]{dghm} and the references 
there.

Consequently, we will only have to prove Theorem~\ref{thm1} for the 
exceptional types $G_2$, $F_4$, $E_6$, $E_7$ and $E_8$. In Section~2,
we present a result (originally due to Geck--Rouquier \cite{gero1}), 
which reduces that proof to the verification of a finite number of cases. 
These finitely many cases will be settled by some explicit computations in 
Section~3, using the {\sf CHEVIE} package \cite{chevie} of the computer 
algebra system {\sf GAP} \cite{gap}. The numbers $|\Irr(H_{\C}(W,\zeta_e))|$ 
(for $W$ of exceptional type) are printed in the table in 
\cite[11.5.13]{ourbuch}.

\begin{rema} \label{rem0} {\rm The conditions on $\ell$ in Theorem~\ref{thm1} 
mean that $\ell$ is a ``bad prime'' for $W$ in the sense of \cite[p.~28]{Ca}.
Simple examples show that, in general, we will have $|\Irr(H_k(W,q))|<
|\Irr(H_{\C}(W,\zeta_e))|$ if the characteristic of $k$ is bad for $W$; 
see, for example, type $G_2$ in \cite[Ex.~3.11]{dghm} and the remarks at
the end of this paper.}
\end{rema}

\section{The generic algebra and its specializations}
We keep the setting and the notation of Section~1. In order to get hold
of $\Irr(H_k(W,q))$, we will work with the generic Iwahori--Hecke algebra 
associated with $W$ and use specialization arguments. For this purpose, we 
have to introduce some notation. We use \cite{ourbuch} as a standard
reference for general facts about finite Weyl groups and Iwahori--Hecke
algebras. 

Let $A={\Z}[v,v^{-1}]$ be the ring of Laurent polynomials in an 
indeterminate~$v$; let $u=v^2$. Then we define the {\it generic 
Iwahori--Hecke algebra} ${\bH}= {\bH}_{A}(W,u)$ in a similar way as 
$H_k(W,q)$ was defined in Section~1. Although now we are not working 
over a field, it is still true that ${\bH}$ is free as an $A$-module, 
with basis $\{T_w \mid w \in W\}$; see \cite[Chapter~4]{ourbuch}. Let $K$ be 
the field of fractions of $A$ and ${\bH}_K=K \otimes_A {\bH}$ be the $K$-algebra 
obtained by extending scalars from $A$ to~$K$; we have ${\bH}_K=H_K(W,u)$
canonically. It is known that the algebra ${\bH}_K$ is split semisimple; 
see \cite[9.3.5]{ourbuch}. 

Let $\tau \colon {\bH} \rightarrow A$ be the $A$-linear map defined by
$\tau(T_1)=1$ and $\tau(T_w)=0$ for $1 \neq w \in W$. Then, by 
\cite[8.1.1]{ourbuch}, $\tau$ is a symmetrizing trace on ${\bH}$ and 
${\bH}$ is a symmetric algebra. Let $\Irr({\bH}_K)$ be the set of simple 
${\bH}_K$-modules, up to isomorphism. Let $c_V \in K$ be the Schur element 
corresponding to $V \in \Irr({\bH}_K)$; see \cite[\S 7.2]{ourbuch}. Since 
${\bH}_K$ is split semisimple, we have $c_V \neq 0$ for all $V$. Since $A$ 
is integrally closed in $K$, we have $c_V \in A$; see 
\cite[7.3.8]{ourbuch}. We shall need the following fact, which is a 
combination of \cite[p.~75]{Ca}, \cite[9.3.6]{ourbuch} and 
the remarks in \cite[3.4]{Lu2}, \cite[2.4]{mykl}.

\begin{abschnitt} \label{schur} {\em Let $P_W=\displaystyle\sum_{w \in W} 
u^{l(w)}$. Then we have a factorization 
\[(u-1)^{|S|}P_W=\prod_{i=1}^{|S|} (u^{d_i}-1), \qquad \mbox{where $d_i
\geq 1$}.\]
The $d_i \geq 1$ are called the degrees of $W$; we have $|W|=d_1 \cdots 
d_{|S|}$. Each Schur element $c_V$ lies in ${\Q}[u,u^{-1}]$ and divides 
$P_W$. The denominators in the coefficients of $c_V$ are divisible by bad 
primes~only.}
\end{abschnitt}

Now let $k$ be a field and $0 \neq q \in k$ be such that $q$ has a square 
root in~$k$. Then there is a unique ring homomorphism $\theta \colon A 
\rightarrow k$ such that $\theta(u)=q$ (and $v$ is mapped to a chosen 
square root of $q$ in $k$). Regarding $k$ as an $A$-module via $\theta$, we 
can extend scalars from $A$ to $k$ and obtain a $k$-algebra ${\bH}_k=k 
\otimes_A {\bH}$ which is canonically isomorphic to $H_k(W,q)$. Thus, it 
will be sufficient to work with ${\bH}$ and its specializations. We can now 
settle a large number of cases occurring in Theorem~\ref{thm1}.

\begin{lem} \label{rem1} Assume that the characteristic of $k$ satisfies
the conditions in Theorem~\ref{thm1}. Then, if $q^{d_i} \neq 1$ for all~$i$
(where the $d_i$ are the degrees of $W$), the algebra ${\bH}_k$ is split 
semisimple and we have $|\Irr({\bH}_k)|=|\Irr({\bH}_K)|=|\Irr(H_{\C}(W,1))|$.
\end{lem}

\begin{proof} Using (\ref{schur}), the above conditions on $k$ and $q$
imply that  $\theta(c_V) \neq 0$ for all $V \in \Irr(\bH_K)$. Consequently,
by \cite[9.3.9]{ourbuch}, the algebra ${\bH}_k$ is split semisimple. Then 
Tits' Deformation Theorem (see \cite[7.4.6]{ourbuch}) applies and we have 
$|\Irr({\bH}_k)|=|\Irr({\bH}_K)|$. Applying the same argument to the 
specialization $A \rightarrow \C$, $u \mapsto 1$, we obtain 
$|\Irr({\bH}_K)|= |\Irr({\bH}_{\C}(W,1))|$. 
\end{proof}

Thus, from now on, we assume that the conditions in Theorem~\ref{thm1}
on the characteristic of $k$ are satisfied and that $q$ is a root of unity. 
Let $k_0 \subseteq k$ be the field of fractions of the image of $\theta$. 
Then, by \cite[3.6]{mykl}, the algebra ${\bH}_{k_0}$ is split. So the scalar 
extension from $k_0$ to $k$ defines a bijection between $\Irr({\bH}_{k_0})$ and 
$\Irr({\bH}_k)$; see \cite[Lemma~7.3.4]{ourbuch}. Thus, we may assume without 
loss of generality that $k=k_0$. 

It will be further convenient to take the following point of view. Let 
$\fp$ be the kernel of $\theta$. Then $\fp \subset A$ is a prime ideal and, 
identifying $k$ with the field of fractions of $A/\fp$, we may regard 
$\theta$ as the natural map $A \rightarrow A/\fp \subseteq k$. Let $e$ be 
the smallest $i \geq 2$ such that $1+q+q^2+\cdots +q^{i-1}=0$ and $\Phi_e(u) 
\in {\Z}[u]$ be the $e$th cyclotomic polynomial. Then we have $\Phi_e(q)=0$ 
in $k$ and so $\Phi_e(u)=\Phi_{e}(v^2) \in \fp$. Now we have 
\begin{alignat*}{2}
\Phi_e(v^2) &=\Phi_{2e}(v)  &&\qquad \mbox{if $e$ is even},\\
\Phi_e(v^2) &=\Phi_e(v)\Phi_{2e}(v) &&\qquad \mbox{if $e$ is
odd}.
\end{alignat*}
Thus, choosing a suitable square root of $q$ in $k$, we can assume that 
$\Phi_{2e}(v) \in \fp$. Let $\fq \subset A$ be the prime ideal generated
by $\Phi_{2e}(v)$; we have $0 \neq \fq \subseteq \fp \subset A$ and
$A/\fq={\Z}[\zeta_{2e}]$, the ring of algebraic integers in the field 
$F={\Q}[\zeta_{2e}] \subset \C$; see \cite[(4.5)]{CR2}. Then, as above, 
we see that ${\bH}_F=H_F(W,\zeta_e)$ is split and the scalar extension 
from $F$ to $\C$ defines a bijection $\Irr({\bH}_F) \leftrightarrow 
\Irr(H_{\C}(W,q))$. 

\begin{abschnitt} {\bf Factorization of decomposition maps.} \label{rem2} 
In the above set-up, the natural map $A \rightarrow A/\fp \subseteq k$ 
induces a well-defined decomposition map $d_{\fp} \colon R_0({\bH}_K) 
\rightarrow R_0({\bH}_k)$ between the Grothendieck groups of ${\bH}_K$ and 
${\bH}_k$; see \cite[7.4.3]{ourbuch}. Similarly, the map $A \rightarrow 
A/\fq \subseteq F$ induces a decomposition map $d_e\colon R_0({\bH}_K) 
\rightarrow R_0({\bH}_F)$. Since $A/\fq$ is integrally closed in $F$,
we have the following factorization (see \cite[2.6]{mybourb}):
\[\mbox{\begin{picture}(160,55)
\put(20,45){$R_0({\bH}_K)$}
\put(65,47){\vector(1,0){80}}
\put(105,52){$\scriptstyle{d_\fp}$}
\put(155,45){$R_0({\bH}_k)$}
\put(48,37){\vector(3,-2){35}}
\put(55,20){$\scriptstyle{d_e}$}
\put(132,13){\vector(3,2){35}}
\put(160,19){$\scriptstyle{d_{\fp}^e}$}
\put(88,05){$R_0({\bH}_F)$}
\end{picture}}\]
Here, $d_{\fp}^e$ is the decomposition map induced by the natural map 
$A/\fq \rightarrow A/\fp$. Now, since $d_{\fp}$ and $d_e$ are surjective 
by \cite[3.3]{mykl}, we conclude that $|\Irr({\bH}_k)|$ equals the rank 
of~$d_{\fp}$ and the above factorization implies that 
\[ |\Irr({\bH}_k)| \leq |\Irr({\bH}_F)|=|\Irr(H_{\C}(W,\zeta_e))|.\]
Note that, if $\fq=\fp$, then Theorem~\ref{thm1} is now clear since 
${\bH}_k={\bH}_F$. 
\end{abschnitt}

Thus, we are finally reduced to the situation where $\fq \neq \fq$. In this 
case, $\fp$ is a maximal ideal containing $\Phi_{2e}(v)$ and a prime number 
$\ell>0$. (See the description of all prime ideals of $A$ in 
\cite[Exercise~7.9]{ourbuch}, for example.)

In order to proceed, we shall need some facts about the center of
${\bH}$, which we denote by $Z({\bH})$. Let $\Cl(W)$ be the set of conjugacy 
classes of $W$. For each $C \in \Cl(W)$, consider the element 
\[z_C=\sum_{w \in W} u^{-l(w)} f_{w,C} T_w \in {\bH},\]
where $f_{w,C}(u) \in {\Z}[u]$ are the class polynomials defined in 
\cite[\S 8.2]{ourbuch}. By \cite[8.2.4]{ourbuch} (see also \cite{gero1}),
we have $z_C \in Z({\bH})$ and $\{z_C \mid C \in \Cl(W)\}$ is an $A$-basis 
of $Z({\bH})$.

\begin{defn}[Geck--Rouquier] \label{def1} {\rm Consider the natural map
$A \rightarrow A/\fq \subseteq F$ and the corresponding decomposition
map $d_e\colon R_0({\bH}_K) \rightarrow R_0({\bH}_F)$. For $M \in 
\Irr({\bH}_F)$, we define $P_{M} \colon Z({\bH}_K) \rightarrow K$ by
\[P_M(z_C)=\sum_{V \in \Irr({\bH}_K)} \frac{d_{V,M}}{c_V} \, \omega_V(z_C)
\qquad \mbox{for $C \in \Cl(W)$},\]
where $d_{V,M}$ are the decomposition numbers and $\omega_V$ is the central
character associated with $V$. By \cite[7.5.3]{ourbuch}, we have 
\[P_M(z_C) \in A_{\fq} \qquad \mbox{for all $C \in \Cl(W)$}.\]
Note that, since $\fq$ is a principal ideal, the localization $A_{\fq} 
\subset K$ is a discrete valuation ring, with residue field~$F$.

For $C \in \Cl(W)$, let $\overline{z}_C:=\sum_{w \in W} \zeta_e^{-l(w)}
f_{w,C}(\zeta_{e}) T_w \in {\bH}_F$. Then $\overline{z}_C \in Z({\bH}_F)$ and 
$\{\overline{z}_C \mid C \in \Cl(W)\}$ is an $F$-basis of $Z({\bH}_F)$; 
see \cite[8.2.5]{ourbuch}. Thus, we obtain an induced function}
\[ \overline{P}_M \colon Z({\bH}_F) \rightarrow F, \qquad \overline{P}_M
(\overline{z}_C)= P_M(z_C) \bmod \fq.\]
\end{defn}

In a slightly different form, the following results appeared in 
\cite[3.3]{gero1}.

\begin{lem} \label{lem1} Assume that each Schur element $c_V$ can be 
expressed as 
\begin{equation*}
c_V=\Phi_{2e}(v)^{d_V} f_V(v) \quad \mbox{where $d_V \geq 0$ and $f_V(v) 
\in A \setminus \fp$}.\tag{$*$}
\end{equation*}
Then we have $P_M(z_C) \in A_{\fp}$ for all $M \in \Irr({\bH}_F)$ and 
$C \in \Cl(W)$. 
\end{lem}

\begin{proof} Let $M \in \Irr({\bH}_F)$. Consider the ideal $I$ of all $g(v) \in 
{\Q}[v,v^{-1}]$ such that $g(v)P_M(z_C) \in {\Q}[v,v^{-1}]$ for all $C  \in
\Cl(W)$. Since ${\Q}[v,v^{-1}]$ is a principal ideal domain, $I$ is 
generated by a polynomial $g_0(v) \in {\Q}[v,v^{-1}]$. Now, since 
$P_M(z_C) \in A_{\fq}$ for all $C \in \Cl(W)$, there exists a polynomial
$h(v) \in I$ which is not divisible by $\Phi_{2e}(v)$. Thus, $\Phi_{2e}(v)$
does not divide $g_0(v)$.

On the other hand, Definition~\ref{def1} shows that the 
product of all $c_V$ lies in~$I$ (note that $\omega_V(z_C) \in A$
since $A$ is integrally closed in~$K$). So, using ($*$) and setting
$d:=\sum_V d_V \geq 0$, we find that
\[ \Phi_{2e}(v)^d f(v) \in I \qquad \mbox{where } f(v)=\prod_{V \in 
\Irr({\bH}_K)} f_V(v) \in A \setminus \fp.\]
Since ${\Q}[v,v^{-1}]$ is a factorial ring, the fact that 
$\Phi_{2e}(v)$ does not divide $g_0(v)$ implies that $g_0(v)$ divides
$f(v)$. Consequently, we have $f(v)P_M(z_C) \in {\Q}[v,v^{-1}]$ for all 
$C \in \Cl(W)$. 

Now (\ref{def1}) actually shows that $\Phi_{2e}(v)^d f(v)P_M(z_C) \in A$ 
for all $C \in \Cl(W)$. Therefore, since $f(v) \in A$ and $\Phi_{2e}(v)$ 
is monic, we conclude that 
\[f(v)P_M(z_C) \in A \qquad \mbox{for all $C \in \Cl(W)$}.\]
Thus, since $f(v) \in A \setminus \fp$, we have $P_M(z_C) \in A_{\fp}$ for
all $C \in \Cl(W)$. \end{proof}

\begin{prop}[Geck--Rouquier] \label{cor1} In the above set-up, assume
that the condition~($*$) in Lemma~\ref{lem1} holds. Then we have 
$|\Irr({\bH}_F)|=|\Irr({\bH}_k)|$.
\end{prop}

\begin{proof} Consider the natural map $A_{\fp} \rightarrow k$ and the 
corresponding decomposition map $d_{\fp} \colon R_0({\bH}_K) \rightarrow 
R_0({\bH}_k)$. For $V \in \Irr({\bH}_K)$, let $\chi_V\colon {\bH}_K 
\rightarrow K$ be the character afforded by $V$. Since $A$ is integrally 
closed in $K$, we have $\chi_V(T_w) \in A$ for all $w \in W$. Thus, we obtain 
an induced function $\tilde{\chi}_V \colon {\bH}_k \rightarrow k$ such that 
$\tilde{\chi}_V(T_w)=\chi_V(T_w) \bmod \fp$. By the definition of $d_{\fp}$, 
the function $\tilde{\chi}_V$ is a $k$-linear combination  of the characters 
afforded by the simple ${\bH}_k$-modules. Hence, we certainly have 
\[ |\Irr({\bH}_k)| \geq \dim_k\, \langle \tilde{\chi}_V \mid V \in \Irr({\bH}_k)
\rangle_k.\]
On the other hand, consider the decomposition map $d_{e}$ induced by
$A_{\fq} \rightarrow F$. Let $\cB \subseteq \Irr({\bH}_K)$ be such that the 
matrix of decomposition numbers $d_{V,M}$, where $V \in \cB$ and $M \in 
\Irr({\bH}_F)$, is square and invertible over $\Z$. Such a subset $\cB$ 
exists by \cite[Theorem~3.3]{mykl}. We have already seen in~(\ref{rem2}) 
that $|\cB|=|\Irr({\bH}_F)|\geq |\Irr({\bH}_k)|$. Hence it is enough to 
show that 
\begin{equation*}
\{\tilde{\chi}_V \mid V \in \cB\} \mbox{ is linearly independent}.
\tag{$\dagger$}
\end{equation*}
To prove ($\dagger$), assume that we have a linear relation
$\sum_{V \in \cB} \lambda_V\, \tilde{\chi}_V=0$ with $\lambda_V \in k$.

Now we consider once more the function $P_M \colon Z({\bH}_K) \rightarrow K$ 
introduced in Definition~\ref{def1}. By Lemma~\ref{lem1}, we have $P_M(z_C) 
\in A_{\fp}$  for all $C \in \Cl(W)$. To simplify notation, for $h=\sum_{w 
\in W} a_w(v)T_w \in {\bH}$ with $a_w(v) \in A$, we write $\tilde{h}:=
\sum_{w \in W} a_w(q^{1/2}) T_w \in {\bH}_k$, where $q^{1/2}$ is the chosen 
square root of $q$. Then, again by \cite[8.2.5]{ourbuch}, we have 
$\tilde{z}_C \in Z({\bH}_k)$ and $\{\tilde{z}_C \mid C \in \Cl(W)\}$ is a 
$k$-basis of $Z({\bH}_k)$. Thus, we obtain an induced function 
\[ \tilde{P}_M \colon Z({\bH}_k) \rightarrow k, \qquad \tilde{P}_M
(\tilde{z}_C)=P_M(z_C) \bmod \fp.\]
For each $V \in \Irr({\bH}_K)$, let $\chi_V^*$ be the unique element in 
$Z({\bH}_K)$ such that $\omega_V(\chi_{V'}^*)=c_V\delta_{V,V'}$ for all 
$V' \in \Irr({\bH}_K)$. Such an element exists by 
\cite[7.2.6]{ourbuch}; explicitly, we have 
$\chi_V^*=\sum_{w \in W} u^{-l(w)}\chi_V(T_w)T_{w^{-1}}$. Since $A$ is 
integrally closed in $K$, we have $\chi_V^* \in Z({\bH})$. Thus, we obtain 
a corresponding element $\tilde{\chi}_V^* \in Z({\bH}_k)$. The defining 
equation in (\ref{def1}) now yields that 
\[ P_M(\chi_V^*)=\sum_{V'\in \Irr({\bH}_K)} \frac{d_{V',M}}{c_{V'}}\,
\omega_{V'}(\chi_V^*)=d_{V,M} \qquad \mbox{for all $M \in \Irr({\bH}_F)$}.\]
Applying this to the above linear relation, we obtain
\[ 0=\sum_{V \in \cB} \lambda_V \tilde{P}_M(\tilde{\chi}_V^*)=
\sum_{V \in \cB} \lambda_V \tilde{d}_{V,M}, \quad \mbox{where 
$\tilde{d}_{V,M}=d_{V,M} \bmod \fp$}.\]
But $\cB$ was chosen such that the matrix $(d_{V,M})$ is invertible over $\Z$. 
Hence $(\tilde{d}_{V,M})$ is invertible over $k$ and so $\lambda_V=0$ for 
all $V \in \cB$, proving $(\dagger$).
\end{proof}

\section{Proof of Theorem 1.1}
By the discussion in Section~2, in order to complete the proof of 
Theorem~\ref{thm1}, it remains to consider the following situation.
The field $k$ is the field of fractions of $A/\fp$, where $\fp$ is a 
maximal ideal containing a prime number $\ell>0$ which is not bad for~$W$.
In particular, $k$ is a finite field of characteristic~$\ell$. Furthermore,
$\fp$ contains the cyclotomic polynomial $\Phi_{2e}(v) \in {\Z}[v]$, where 
$e \geq 2$ divides some $d_i$ as in (\ref{schur}). The number~$e$ is 
minimal such that $1+q+q^2+\cdots +q^{e-1}=0$ in~$k$ and there is a square 
root of~$q$ which is also a root of $\Phi_{2e}(v)$ in~$k$. Under these 
conditions, we must show:
\begin{equation*}
|\Irr({\bH}_k)|=|\Irr({\bH}_F(W,\zeta_e))|, \qquad \mbox{ where 
$F={\Q}[\zeta_{2e}]\subset {\C}$}.  \tag{E}
\end{equation*}
We shall need the following well-known property of cyclotomic polynomials.

\begin{lem} \label{cyc} Assume that $\Phi_d(q)=0$ for some $d \geq 2$.
If $q=1$, then $d=\ell^n$ for some $n \geq 1$. If $q \neq 1$, then $e$
is the order of $q$ in $k^\times$ and $d=e\ell^n$ for some $n \geq 0$.
\end{lem}

\begin{proof}  Let us write $d=d'\ell^n$ where $\ell$ does not divide $d'$.
Now $\Phi_{d}(u)$ certainly divides $\Phi_{d'}(u^{\ell^n})$. Since, 
furthermore, taking $\ell$th powers is an automorphism of $k$, we conclude 
that $\Phi_{d'}(q)=0$. If $d'=1$, this means that $q=1$ and $e=\ell$.
So the assertion is true in this case. Now assume that $q \neq 1$. 
Then $e$ is minimal such that $q^e-1=0$, i.e., $e$ is the multiplicative 
order of~$q$. Hence $e$ divides the order of $k^\times$ and so is coprime 
to~$\ell$. Then $\Phi_{d'}(q)=0$ implies $q^{d'}=1$ and so $e$ divides~$d'$. 
Assume, if possible, that $e \neq d'$. Then we can write
\[ \frac{u^{d'}-1}{u-1}=a(u)\Phi_{e}(u)\Phi_{d'}(u) \qquad \mbox{where 
$a(u) \in {\Z}[u]$}.\]
Differentiating with respect to $u$ and applying the natural map
$A \rightarrow A/\fp \subseteq k$ yields that $d'=0$ in $k$, contradicting 
the fact that $d'$ is coprime to~$\ell$. \end{proof}

The following result (first established in \cite{gero1}) provides a general
proof of (E), but under a condition on~$\ell$ which is more restrictive 
than that in Theorem~\ref{thm1}.

\begin{thm}[Geck--Rouquier] \label{lem3} Assume that $\ell \neq 2$ and that
$e\ell$ does not divide any degree $d_i$ of $W$. (This condition is 
satisfied, for example, if $\ell$ does not divide $|W|$.) Then condition 
($*$) in Lemma~\ref{lem1} is satisfied and (E) holds.
\end{thm}

\begin{proof} By Proposition~\ref{cor1}, it is enough to verify that 
condition ($*$) in Lemma~\ref{lem1} is satisfied. Using (\ref{schur}), 
we can write $c_V=\Phi_{2e}(v)^{d_V}f_V(v)$, where $d_V \geq 0$ and $f_V(v) 
\in {\Z}[v,v^{-1}]$ is a product of an integer divisible by bad primes 
only, an integral power of~$u$, $\Phi_e(v)$ (if $e$ is odd) and various 
cyclotomic polynomials $\Phi_{d}(u)$, where $d \neq e$ and $d \geq 2$ 
divides some degree~$d_i$.

Now consider the natural map $\theta \colon A \rightarrow A/\fp \subseteq k$.
Assume, if possible, that $\theta(f_V(u))=0$ for some $V \in \Irr({\bH}_K)$. 
Since the characteristic of $k$ is not a bad prime, we must have 
$\theta(\Phi_e(v))=0$ (if $e$ is odd) or $\theta(\Phi_d(u))=0$ for some 
$d \neq e$ where $d \geq 2$ divides a degree~$d_i$. The second possibility
cannot occur by Lemma~\ref{cyc} and our assumption. Hence $e$ must be odd 
and $\theta(\Phi_e(v))=0$. As in the proof of Lemma~\ref{cyc}, we write
$(v^{2e}-1)/(v-1)=b(v)\Phi_{e}(v)\Phi_{2e}(v)$ with some $b(v) \in {\Z}[v]$,
differentiate with respect to $v$ and obtain the conclusion that 
$\theta(2e)=0$, i.e., $\ell$ divides $2e$. This contradicts the fact
that $\ell \neq 2$ and $e$ is the order of $q$. 
\end{proof}

In Section~1, we have already remarked that Theorem~\ref{thm1} is known
to hold for $W$ of type $A_n$, $B_n$ and $D_n$. Thus, it remains to consider
the finitely many exceptional types. For each of these types, we need to
know the prime divisors of $|W|$, the bad primes and the degrees $d_i$.
This information is provided in Table~2. (The information in this table 
is obtained by inspection of the tables in \cite[Appendix~E]{ourbuch}
and \cite[p.~75]{Ca}.)

\begin{table}[htbp]  \caption{Bad primes and degrees for the exceptional 
types.}
\begin{center}
$\begin{array}{cccc} \hline 
\mbox{Type} & |W| & \mbox{$\ell$ bad} & \mbox{degrees $d_i$}\\ \hline
G_2 & 2^2 \cdot 3 & 2,3 & 2,6 \\
F_4 & 2^7 \cdot 3^2 & 2,3 & 2,6,8,12 \\
E_6 & 2^7 \cdot 3^4 \cdot 5 & 2,3 & 2,5,6,8,9,12 \\  
E_7 & 2^{10}\cdot 3^4 \cdot 5 \cdot 7 & 2,3 & 2,6,8,10,12,14,18\\
E_8 & 2^{14} \cdot 3^5 \cdot 5^2 \cdot 7 & 2,3,5 
& 2,8,12,14,18,20,24,30 \\ \hline
\end{array}$
\end{center}
\end{table}

Now, given $W$, $e$ and $\ell$, the {\sf CHEVIE} function in Table~3 
computes the number of simple $\bH_k$-modules. This is 
done by computing the rank of the character table of $\bH_K$ under the 
specialization $A \rightarrow A/\fp \subseteq k$. That this rank indeed 
equals the cardinality of $|\Irr(\bH_k)|$ follows from \cite[7.5.7 and 
8.2.9]{ourbuch} and the fact that $d_{\fp}$ is surjective (see 
\cite[3.3]{mykl}). Note also that the character tables for all $\bH_K$
are explicitly known and available in {\sf CHEVIE}.

By Theorem~\ref{lem3}, we still have to consider all primes $\ell$ which 
divide $|W|$ but are not bad and all $e$ such that $e$ and $e\ell$ divide 
some degree $d_i$. Thus, if $W$ is of type $G_2$ or $F_4$, then all prime 
divisors of $|W|$ are bad and there is nothing more to prove. If $W$ is of 
type $E_6$, we only have to consider the prime $\ell=5$. But the only
degree divisible by $5$ is $5$ itself, hence we are done. If $W$ is of 
type $E_7$ or $E_8$, then $2$ is the only possible value for $e$ such 
that $e$ and $e\ell$ (with $\ell$ bad) divide some $d_i$. Hence it 
remains to consider the following cases:
\begin{itemize}
\item[$E_7$:] $e=2$ and $\ell=5,7$, where we obtain $|\Irr(\bH_k)|=12$.
\item[$E_8$:] $e=2$ and $\ell=7$, where we obtain $\Irr(\bH_k)|=23$.
\end{itemize}
These results are in accordance with the corresponding values in
characteristic~$0$ printed in \cite[11.5.13]{ourbuch}. Thus, 
Theorem~\ref{thm1} is proved.

\begin{table}[htbp] \caption{A {\sf CHEVIE} program for computing the number
of simple modules.}
{\small \begin{verbatim} 
  gap>  RequirePackage("chevie");
  gap>  NumberSimples := function( W, e, ell )
  gap>     local H, i, z;
  gap>     if  Gcd( ell, e ) = 1  then
  gap>       # find i such that GF(ell^i) has an element z 
  gap>       #of order 2e
  gap>       i:=1;
  gap>       while  not IsInt( ( ell^i-1 )/( 2*e ))  do
  gap>          i := i+1;
  gap>       od;
  gap>       z := Z( ell^i )^( ( ell^i-1 )/( 2*e ) );
  gap>     elif  e = ell  and  ell>2  then
  gap>       z := Z( ell )^( ( ell-1 )/2 );
  gap>     elif  e = ell  and  ell=2  then
  gap>       z := Z( ell );
  gap>     fi;
  gap>     # define the Iwahori-Hecke algebra with parameter z^2
  gap>     H := Hecke( W, z^2, z );
  gap>     # return the rank of the specialized character table
  gap>     return RankMat( CharTable( H ).irreducibles );
  gap>  end;
  gap>  NumberSimples( CoxeterGroup( "E", 8 ), 2, 5 );  # an example
  22
\end{verbatim}}
\end{table}

\begin{rema} {\rm Using the above program we find that if $e=\ell$ is 
a bad prime, then $|\Irr(\bH_k)|< |\Irr(\bH_F)|$. So, in general, the 
assertion in Theorem~\ref{thm1} will not hold if $\ell$ is a bad prime. 
(Another example is given in Table~3.)}
\end{rema}

{\sc\footnotesize Institut Girard Desargues, bat. 101, Universit\'e 
Claude Bernard Lyon 1, 43 bd du 11 novembre 1918, F--69622 Villeurbanne cedex, 
France}

\smallskip
\makeatletter
{\footnotesize\it E-mail address: \tt geck@desargues.univ-lyon1.fr}
\makeatother

\end{document}